\documentclass[11pt]{amsart}
\usepackage{graphicx}
\usepackage{url} 
\usepackage{hyperref} 
\newtheorem{thm}{Theorem}[section]
\newtheorem{cor}[thm]{Corollary}
\newtheorem{lem}[thm]{Lemma}
\newtheorem{prop}[thm]{Proposition}
\theoremstyle{definition}
\newtheorem{defn}[thm]{Definition}
\theoremstyle{remark}

\numberwithin{equation}{section}
\begin{document}

\title[Unbounded Collatz Sequence]{On the Behavior of Unbounded Collatz Sequences}%
\author{ Jorge Salazar}

\address{DMAT, Universidade de \'Evora, \'Evora - Portugal}%
\email{salazar@uevora.pt}%

%\thanks{}
\subjclass{Number Theory}%

% ----------------------------------------------------------------
\begin{abstract}
 The aim of this paper is to show a peculiar behavior  of a (hypothetical)  Collatz sequence  going to infinity. We study the associated  Syracusa sequence (the odd elements of the former) and show that the limit set of a conveniently  normalized sequence is the whole unit interval. In particular, for  any positive integer there is a subsequence whose elements' expansions in base 3  begin (from the left) with the expansion of the given number.

\end{abstract}
\maketitle
% ----------------------------------------------------------------

% ----------------------------------------------------------------
\section{Introduction}

 \noindent
 The Collatz operation is $Col(a)=3a+1$,  if $a$ is odd and $Col(a)=a/2$ if $a$ is even. 
% The Collatz process is the sequence of the iterations of this simple operation.
 In 1937  Lothar Collatz introduced this operation and conjectured that any sequence of its iterations eventually becomes the loop $ 1,4,2,1 $. The  conjecture  is also known as the “3x+1 conjecture” or the “Syracuse problem.” Terence Tao \cite{Taopres} gave an enthusiastic slide presentation of this conjecture and he himself have important contributions \cite{Tao1}.
 
  \vskip 5pt
 \noindent
Since  $3a+1$ is  even, we can divide by $2$ in the same step and consider instead the iterations of
\[ Col_2(a)=\left\lbrace
\begin{array}{cl}
\frac{1}{2}\left( 3a+1\right) ,& \hbox{if } a \hbox{ is odd},\\
&\\
\frac{1}{2}a,& \hbox{if } a \hbox{ is even}.
\end{array} \right.  \]

  \vskip 5pt
 \noindent
Given  $a_0\ge 1$, denote $  a_{\alpha}=Col_2\left( a_{\alpha-1}\right) ,\ \alpha\ge 1 $,
%	\begin{equation}\label{aa}\end{equation}
 the sequence of iterations. Let's denote $q_\alpha$ the biggest power in the development of $a_\alpha$ in base 3. i.e.
\begin{equation}\label{qalfa}
1\le 3^{-q_\alpha}a_\alpha < 3 % =\left[ \frac{\ln a_\alpha}{\ln 3} \right] ,
\end{equation}
% where $\left[ \, \cdot\,\right] $ denotes the biggest integer less than or equal to the given  number.
 
%  \vskip 5pt \noindent 
% The Collatz conjecture says that any orbit contains 1, eventually becoming the loop 1,2,1. The book  by Jeffrey C. Lagarias (ed.) \cite{laga} is an exelent source, we can appreciate so many different approaches, from Number Theory, Dynamical Systems, Stochastic Processes,  Numerical Methods, also including the motivation and origin of this problem as well as   an annotated bibliography from 1963 to1999. The very recent article by Terence Tao \cite{Tao}, refined considerably  previous results on this problem, from the probabilistic perspective and revived the  interest of the public.  Following a seemingly overlooked  approach, the aim of this paper is to prove the theorem below. 
 
 {\thm \label{collatz} Assume that for some integer $a_0\ge 1$, the Collatz sequence verifies $\lim_{\alpha \rightarrow \infty} a_\alpha=\infty$. Let $\psi $ be any positive integer. Then there are infinitely many indexes $\alpha$ such that $ \left[ 3^{p-q_\alpha} a_\alpha \right] =\psi$, where  $p=\left[ \frac{\ln \psi}{\ln 3} \right]$.} 
 
 \vskip 5pt
 \noindent 
 ($\left[ \cdot \right] $ is the integer part of a given number.)

 \vskip 5pt
\noindent 
The equality  $\left[ 3^{p-q_\alpha} a_\alpha \right] =\psi$  means that the first $p$ digits (from the left) of  $a_\alpha$ in base $3$ match those of $\psi$.

 \vskip 5pt
\noindent 
In the whole paper, $\left( a_\alpha \right)_{\alpha\in \mathbb{{N}}} $ is a \textbf{fixed} orbit of $Col_2$ going to infinity. i.e.
\begin{equation}\label{key10}
\forall \alpha\ge 1, \,  a_\alpha =Col_2\left(  a_{\alpha -1}\right) \hbox{ and } \lim_{\alpha \rightarrow \infty} a_\alpha=\infty.
\end{equation}

\section{Bounded  version of the Collatz process}

\noindent
For all $\alpha\ge 1$, put
\begin{equation}\label{c}
c_\alpha=3^{-q_\alpha}\cdot a_\alpha .
\end{equation}
The following definitions are intended to obtain  $c_\alpha $ as the iterations of a single operation starting from $c_1$. i.e. We need to understand how the Collatz operation works on $c_\alpha$, for $\alpha\ge 1$.

\vskip 5pt
\noindent
\defn  Denote
\begin{equation}\label{key11}
\mathbf{N} =\bigcup_{r=0}^{\infty}3^{-r} \mathbb{N}.
\end{equation}
the set of all the rational numbers with a finite expansion  in base 3

\vskip 5pt
\noindent
For $b\in \mathbf{N}$, define
 \begin{equation}\label{ctilde1}
Col_3(b)=\left\lbrace 
\begin{array}{ll}
\frac{b}{2}, & \hbox{if } \frac{b}{2}\in \mathbf{N},\\ &\\
\frac{b}{2} +  \frac{3^{-q-1}}{2}  ,& \hbox{if } \frac{b}{2}\notin \mathbf{N},
\end{array}\right. 
\end{equation}
where 
\begin{equation}\label{key12}
q=\min\left\lbrace r \ge 0;\, 3^rb\in \mathbb{N}\right\rbrace .
\end{equation}

\vskip 5pt
\noindent 
\emph{Note}: If  $\frac{b}{2}\notin \mathbf{N}$, the base 3 expansion of $\frac{b}{2}$ ends in an infinite string of $1$'s, starting at the position  ${-q-1}$.
 In that case,  we ``round up" to $\frac{b+3^{-q-1}}{2}  $  and   $ Col_3(b)\in \mathbf{N} $.  When the rounding occurs, the last position is ${-q-1}$ and its value is $2$.
 Working with $Col_3$, we take advantage of the long division procedure (in base 3) by starting with division by 2, without caring about the parity of the number until we get to the last digit.
 
 \vskip 5pt
 \noindent 
For future reference,  observe that for all ${b}\in \mathbf{N}$ and $n\in\mathbb{N}$,
 \begin{equation}\label{key200}
 {2}^{-n} b \le Col_3^{n}\left( b\right) \le    {2}^{-n} b + 3^{-q }   , %\cdot {2}^{-n_j},
 \end{equation}
 where $q$ is defined by \eqref{key12}.

\vskip 5pt
 \noindent
\defn  For all  $b\in \mathbf{N}\cap \left[ 1,3\right) $, define
 \begin{equation}\label{ctilde}
 Col_4(b)=\left\lbrace \begin{array}{rl}
 Col_3(b) , & \hbox{if } 2\le b <3 ,\\ &\\
 3\cdot Col_3(b) , & \hbox{if }  1\le  b< 2 .
 \end{array} \right. 
 \end{equation}

  \vskip 5pt
 \noindent
 Clearly, for all $\alpha \ge 1$,
 \begin{equation}\label{cc}
 c_{\alpha +1}=Col_4\left(c_\alpha \right) .
 \end{equation}
(If $a_0$ is not divisible by $3$, \eqref{cc} is also true for $\alpha =0$.)
 
 \vskip 5pt
 \noindent
 The usefulness of this version of the Collatz iteration ($Col_4$) is  that we can multiply by  the appropriate power of $3$ after any number of iterations of $Col_3$, as we can see in the next Lemma. 
 
 {\lem\label{col4}  For all $n\in\mathbb{N}$, for all  $b\in \mathbf{N}\cap \left[ 1,3\right) $, 
 \[ Col_4^n(b)=3^k\cdot Col_3^n(b) , \]
 where 
   \begin{equation}\label{nk1}
 k=-\left[ \frac{\ln  Col_3^n(b)}{\ln 3} \right] .
 \end{equation}}
%\[   1\le 3^k\cdot\widetilde{C}^n(b) <3.\]

 \vskip 5pt
 \noindent
Indeed, for all $\alpha$, the coeficients of  $a_\alpha$, $c_\alpha$  and  
\[  b_{\alpha+1}=Col_3\left( b_{\alpha}\right) ,\ \alpha\ge 1, \ b_1=c_1,  \]
in base $3$, are the same, just shifted. % (assuming $a_0$ is not divisible by $3$).

 \vskip 5pt
\noindent 
The main ingredient in the proof of Theorem  \ref{collatz} is the following result (interesting on their own), which will be proved in Section \ref{proof}.

{\prop \label{c=1}
	\begin{equation}\label{ccc}
	\liminf_{\alpha\rightarrow\infty}c_\alpha=1.
	\end{equation}}

\section{preliminary results}

{\lem\label{nk2} Let $  \phi , \theta  \in  \mathbf{N}\cap \left[ 1,3\right) $. Suppose there is an integer $k $ such that 
	\begin{equation}\label{phi}
	\left[ 3^k\phi\right] =\left[ 3^k \theta\right] . 
	\end{equation}
Assume that $\left[ 3^k\phi\right] \neq3^k\phi$ and $ \left[ 3^k \theta\right] \neq  3^k \theta$ (irrelevant if  $\phi=\theta$). Then, for all $m\in \mathbb{N}$,
\begin{equation}\label{longdiv}
	\left[ \frac{3^k}{2^{n}} \phi\right] = \left[ 3^kCol_3^m\left( \phi\right) \right] =\left[ 3^k Col_3^m\left( \theta\right) \right] =\left[ \frac{3^k}{2^{n}}  \theta \right]  .
	\end{equation} }

\vskip 5pt
\noindent
\proof   Write $  \phi $ and $ \theta $    in base 3,  and divide them by $2$ using the long division procedure. % (i.e. from left to right). 
\vskip 5pt
\noindent
(The  assumption   $\left[ 3^k\phi\right] \neq3^k\phi$ and $ \left[ 3^k \theta\right] \neq  3^k \theta$ ensures that the roundings on all the iterations occur outside of the first $k$ digits (in base 3).  %$0\cdots k$.) 
\vskip 5pt
\noindent
\emph{Observation}: Since $Col_3^m(b)\rightarrow 0 $, when $m\rightarrow\infty$, the  Lemma \ref{nk2} becomes irrelevant for $m$ big enough. To keep it meaningful, we restrict it to $m\le n$, where $n$ is defined  by 
\begin{equation}
\left[ 3^kCol_3^n\left( \phi\right) \right] =\left[ 3^k Col_3^n\left( \theta\right) \right]  =1.
\end{equation}

\vskip 5pt
\noindent
After the next iteration, all the meaningful  digits of  $  Col_3^{n+1}\left( \phi\right)$ and $ Col_3^{n+1}\left( \theta\right)$  (in base 3) disappear from the  ``window" we are observing (the first $k$  positions)   and  there is nothing more we can say from \eqref{phi}.
\[ \left[ 3^kCol_3^{n+1}\left( \phi\right) \right] =\left[ 3^k Col_3^{n+1}\left( \theta\right)\right]  =0. \]

\defn\label{Phi} Given an integer $\xi\ge 1$,  put  $ q_0=\left[ \frac{\ln\xi}{\ln 3}\right] $ and denote  
\begin{equation}\label{key81}
\Phi_{\xi}=\left\lbrace   b\in \mathbf{N}\setminus \bigcup_{r=0}^{q_0} 3^{-r}\mathbb{N} ;\, \left[ 3^{q_0}b\right] =\xi\right\rbrace 
\end{equation}
the set of all the numbers $b\in \left[ 1,3\right) $, whose digital expansion in base 3 matches the digits in the expansion of $\xi$, in the same order, but in a shifted position. In particular,  $\Phi_{\xi}\subset \left[ 1,3\right) $.

{\cor\label{cor1}
For $b\in\Phi_{\xi}$ and $n\in\mathbb{N}$, define
 \begin{equation}\label{key82}
 	k_1=-\left[ \frac{\ln  Col_3^n(b)}{\ln 3} \right]  \hbox{ and }  \, k_2=-\left[ \frac{\ln b}{\ln 3}  - \frac{n\ln2}{\ln 3} \right] .
 \end{equation}
Then, as long as $k_1\le q_0$ or $k_2\le q_0$, we have $k_1=k_2$ and their common value, $k$, does not depend on $b\in\Phi_{\xi}$, only on $n$. 

\vskip 5pt
\noindent
In particular, for all pairs $\left( n,k\right) $, $k\le q_0$ verifying \eqref{key82},  
 \begin{equation}\label{key83}
 Col_4^n(b)=3^k\cdot Col_3^n(b)
 \end{equation}  
for all $b\in\Phi_{\xi}$. }(See Lemma \ref{col4}.)

  \vskip 5pt
 \noindent
 For further reference, let's write this conclusion in other form.

{\cor\label{cor2} 
 The pairs $(n,k)$ in  Corollary \ref{cor1}  are exactly the solutions of 
 \begin{equation}\label{pairs}
  1 \le 3^{k}\cdot Col_3^n\left( b\right) <3  \, \hbox{ and/or }\, 1 \le 3^{k}\cdot b\cdot 2^{-n} <3  .
 \end{equation}
 for which $k\le q_0$.}

 \vskip 5pt
\noindent
 Note that  given $n$, $k$ is unique, while for each $k$, there might be  two values of $n$ verifying \eqref{pairs}, but not more than  two.

 {\lem\label{kjnj}  There is a  sequence $\left( n_j, k_j\right) $  (increasing in both components), such that for all real number $x$, $1<x<3$, there is $j_0\in\mathbb{N}$ such that for all $j\ge j_0$,  
 	\begin{equation}\label{pairs2}
 k_j\frac{\ln3}{\ln2}  < n_j\le k_j\frac{\ln3}{\ln2} +\frac{\ln x}{\ln2} .
 	\end{equation}
(In particular, $\left( n_j, k_j\right) $ verifies $ 1 \le 3^{k}\cdot x\cdot 2^{-n} <3  $, see \eqref{pairs}.)

 \vskip 5pt
 \noindent
 \proof Since $\frac{\ln 3}{\ln2} $ is irrational,
 \[\limsup_{k\rightarrow \infty, \, k\in \mathbb{N}} \left(  k\cdot \frac{\ln 3}{\ln2} - \left[  k\cdot \frac{\ln 3}{\ln2} \right]  \right) =1 .\]
 %( $\left[ z\right] $ denotes the biggest natural number $\le z$.)  
 So, there is  an increasing  sequence $k_j \in \mathbb{N} $ such that 
 \[\lim_{j\rightarrow \infty}  k_j\cdot \frac{\ln 3}{\ln2} - \left[  k_j\cdot \frac{\ln 3}{\ln2} \right]  =1 .\]
 Put  \[ n_j= \left[  k_j\cdot \frac{\ln 3}{\ln2}\right] +1 . \]
So, 
 \begin{equation}\label{key99}
  k_j\cdot \frac{\ln 3}{\ln2} < n_j < k_j\cdot \frac{\ln 3}{\ln2} +1.
 \end{equation}
 Since   $\ln x >0$ and
 \[\lim_{j\rightarrow \infty} n_j-  k_j\cdot \frac{\ln 3}{\ln2} =0 ,\]
 there exist $j_0$ such that for all $j\ge j_0$,
 \begin{equation}\label{up}
  n_j <  k_j\cdot \frac{\ln 3}{\ln2} + \frac{\ln x}{\ln2}  . \end{equation}
 By \eqref{key99}and \eqref{up}, we have \eqref{pairs2}.

\section{Proofs}\label{proof}
  \vskip 5pt
% \proof[\textbf{Proof of Theorem \ref{c=1}}] 
\subsection{Proof of Proposition  \ref{c=1}}

Assume that 
 $ l=\liminf_{\alpha\rightarrow\infty}c_\alpha>1 $.

 \vskip 5pt
 \noindent
 Let $\left( n_j, k_j\right) $  be a  sequence  given by Lemma  \ref{kjnj} and fix  $0<\epsilon< \frac{\ln l}{\ln2} $. By  \eqref{pairs2}, there is $j_0$ such that for all $j\ge j_0$,
 \begin{equation}\label{dense} 
 k_{j}\cdot \frac{\ln 3}{\ln2}  <n_{j} <  k_{j}\cdot \frac{\ln 3}{\ln2} + \epsilon .
 \end{equation}
 
 \vskip 5pt
 \noindent
 For later reference, note that
 \begin{equation}\label{densek1} 
 n_{j} \cdot \frac{\ln2}{\ln 3}- \epsilon\cdot \frac{\ln2}{\ln 3}  <k_{j} <n_{j} \cdot \frac{\ln2}{\ln 3}   .\end{equation}

\vskip 5pt
\noindent
Since $\epsilon\cdot \frac{\ln2}{\ln 3}< \frac{\ln l}{\ln3} <1$,
\begin{equation}\label{densek} 
 k_{j} =\left[ n_{j} \cdot \frac{\ln2}{\ln 3}  \right]  .
\end{equation}

 \vskip 5pt
 \noindent
 Since $ \lim_{j\rightarrow \infty}   3^{-{k_{j}}} \left[3^{k_{j}}l-1 \right] =l $, we choose $j_1\ge j_0$  such that  
 \begin{equation}\label{jota}
 \ln  3^{-{k_{j_1}}}\left[3^{k_{j_1}}l -1\right] > \epsilon \ln2 .
 \end{equation}
 
 \vskip 5pt
 \noindent
  Now, choose a sequence   $\left( c_{\alpha_i}\right) _{i\in\mathbb{N}} $, converging to $l$, such that for all $i$,
   \begin{equation}\label{ele}
  \left[ 3^{k_{j_1}} c_{\alpha_i} \right] =\left[3^{k_{j_1}}l \right] .
  \end{equation}

\vskip 5pt
\noindent 
If no such sequence exist,  necessarily   $l\in\mathbf{N}$ and  any  sequence $ c_{\alpha_i}$ going to $ l$ approach $l$ by below.  In that case,  we can choose $ c_{\alpha_i}$  such that for all $i$,
 \begin{equation}\label{ele1}
 \left[ 3^{k_{j_1}} c_{\alpha_i} \right] =\left[3^{k_{j_1}}l -1\right] .
 \end{equation}
The key point  is we have chosen $ c_{\alpha_i}\rightarrow l$  such that $  \left[ 3^{k_{j_1}} c_{\alpha_i} \right]  $ does not depend on $i$. 
In other words, taking $\xi= \left[ 3^{k_{j_1}} c_{\alpha_i} \right]  $ (for some $i$), 
\begin{equation}\label{key101}
 \left( c_{\alpha_i}\right) _{i\in\mathbb{N}} \subseteq \Phi_{\xi} .
\end{equation}

\vskip 5pt
\noindent
Now, observe that 
\begin{equation}\label{ces}
c_{\alpha_i} \ge 3^{-{k_{j_1}}}\left[3^{k_{j_1}}  c_{\alpha_i} \right] \ge 3^{-{k_{j_1}}}\left[3^{k_{j_1}}l -1\right].
\end{equation}
Then, by \eqref{jota}, for all $i$,
\begin{equation}\label{epsy}
\ln c_{\alpha_i}> \epsilon \ln2 
\end{equation}
and, by \eqref{densek1},
\begin{equation}\label{densek2} 
 -k_{j_1}< \frac{\epsilon\ln2}{\ln 3} - \frac{n_{j_1}\ln2}{\ln 3} < \frac{\ln c_{\alpha_i}}{\ln 3} - \frac{n_{j_1}\ln2}{\ln 3}  .
  \end{equation}

\vskip 5pt
\noindent
Therefore, by Corollary \ref{cor1}, 
	\begin{equation} \label{key84}
	k=-\left[ \frac{\ln  Col_3^{n_{j_1}}\left( c_{\alpha_i}\right) }{\ln 3} \right] =-\left[ \frac{\ln c_{\alpha_i}}{\ln 3}  - \frac{n_{j_1}\ln2}{\ln 3} \right] \le k_{j_{1}}
	\end{equation}
	is constant and, by  \eqref{key83},
		\begin{equation} \label{key86}
	Col_4^{n_{j_1}} \left( c_{\alpha_i}\right) =3^{k}\cdot Col_3^{n_{j_1}}\left( c_{\alpha_i}\right) \le 3^{k_{j_1}}\cdot Col_3^{n_{j_1}}\left( c_{\alpha_i}\right) .
	\end{equation}  

	 \vskip 5pt
	\noindent
By \eqref{key200},
\begin{equation}\label{key20}
{2}^{-n_{j_1}} c_{\alpha_i} \le Col_3^{n_{j_1}}\left( c_{\alpha_i}\right) \le    {2}^{-n_{j_1}} c_{\alpha_i} + 3^{-q_{\alpha_i}}   , %\cdot {2}^{-n_j} ,
\end{equation}
	where $q_{\alpha_i} $ is given by    \eqref{c} (or equivalently by \eqref{key12}). Note that $q_{\alpha_i} \rightarrow \infty$, since $  a_\alpha \rightarrow \infty  $.

	 \vskip 5pt
	\noindent
Putting  \eqref{key86} and \eqref{key20} together,
\begin{equation} \label{fin1}
Col_4^{n_{j_1}} \left( c_{\alpha_i}\right) \le \frac{3^{k_{j_1}}}{ {2}^{n_{j_1}} } c_{\alpha_i} + 3^{k_{j_1}-q_{\alpha_i}}   . %\cdot {2}^{-n_j},
\end{equation}  
 
\vskip 5pt
\noindent 
 Finally,  since
 \begin{equation} \label{key44}
 Col_4^{n_j}\left( c_{\alpha_i}\right) = c_{\alpha_i+n_j} 
 \end{equation}
 then
 \begin{equation} % \label{key5}
 l\le \liminf_{i\rightarrow\infty} c_{\alpha_i+n_j}\le \frac{3^{k_j}}{ {2}^{n_j}} l .
 \end{equation}
 Which is a contradiction, since $\frac{3^{k_j}}{ {2}^{n_j}} <1$.
 
 \vskip 5pt
\subsection{Proof of Theorem  \ref{collatz}}
Given an integer $\psi\ge 1$, set  
\begin{equation}\label{key47}
p=\left[ \frac{\ln \psi}{\ln 3} \right]  \hbox{ and  }\, \psi_0= 3^{-p}\psi.
\end{equation}
Put 
\begin{equation} \label{key66}
\tau= \frac{\ln \psi_0}{\ln2} -\left[ \frac{\ln \psi_0}{\ln2 } \right]\left( \frac{\ln 3}{\ln 2}-1 \right) .
\end{equation}
In other words,
\begin{equation} \label{key66b}
\tau= \left\lbrace 
\begin{array}{ll}
\frac{\ln \psi_0}{\ln2} ,& \hbox{if } \left[ \frac{\ln \psi_0}{\ln2 } \right] =0,
\\ 
&\\
\frac{\ln \psi_0}{\ln2} +1-\frac{\ln 3}{\ln 2} , & \hbox{if } \left[ \frac{\ln \psi_0}{\ln2 } \right] = 1.
\end{array} \right. .
\end{equation}
Since $1\le \psi_0 < 3$, $\tau$ is completely defined and $0\le\tau <1 $.

\vskip 5pt
\noindent
By the irrationality of $\frac{\ln3}{\ln2}$, there is a sequence $k_j$ such that  
\begin{equation}\label{key6}
\lim_{j\rightarrow \infty}  k_j\cdot \frac{\ln 3}{\ln2} - \left[  k_j\cdot \frac{\ln 3}{\ln2} \right]  =  \tau.
\end{equation}

\vskip 5pt
\noindent
 For all $\epsilon>0$, there is $j_0$  such that for all $j\ge j_0$,
\begin{equation}\label{key45}
2^{-\epsilon}\cdot 2^\tau< \frac{3^{k_j}}{2^{n_j }}< 2^\epsilon\cdot 2^\tau ,
\end{equation}
where 
\begin{equation}\label{key68}
n_j= \left[  k_j\cdot \frac{\ln 3}{\ln2}\right]  .
\end{equation}

\vskip 5pt
\noindent
Now, fix 
\begin{equation}\label{epsi}
0<\epsilon<\frac{\ln\left( 1+3^{-p-2}\right) }{\ln2}.
\end{equation}
Since 
\begin{equation} \label{key66c}
2^\tau= \left\lbrace 
\begin{array}{ll}
 \psi_0 & \hbox{if } 1\le \psi_0<2,
\\ 
&\\
\frac{2}{3} \psi_0& \hbox{if } 2\le \psi_0<3.
\end{array} \right. ,
\end{equation}
multiplying  all members of \eqref{key45} by $3^p$, for all $j$ big enough, we obtain
\begin{equation}\label{key46}
\psi -\frac{1}{3}  < 2^{-\epsilon}\psi< \frac{3^{k_j+p}}{2^{n_j }}< 2^\epsilon\psi< \psi +\frac{1}{3}  
\end{equation}
if $ 1\le \psi_0<2 $, or 
\begin{equation}\label{key48}
\psi -\frac{1}{3} < 2^{-\epsilon}\psi< \frac{3^{k_j+1+p}}{2^{n_j+1 }}< 2^\epsilon\psi < \psi +\frac{1}{3}
\end{equation}
otherwise.

\vskip 5pt
\noindent
Now, fix $j$ verifying  \eqref{key45} and  take $\xi=3^{k_{j}+q}$ in Definition \ref{Phi},  $q\in \mathbb{N}$ will be determined later. By Proposition \ref{c=1}, there is a subsequence 
\begin{equation}\label{key50}
\left( c_{\alpha_i}\right) _{i\in \mathbb{N}}\subset \Phi_{3^{k_{j}+q}}.
\end{equation}

\vskip 5pt
\noindent
By Corolary \ref{cor1}, for all $i$,
\begin{equation} \label{key79}
Col_4^{n_j}\left( c_{\alpha_i}\right) = 3^{k_j}Col_3^{n_j}\left( c_{\alpha_i}\right) .
\end{equation}
Since $ 1< c_{\alpha_i} < 1 + 3^{-k_j-q}$,
\begin{equation}\label{key21}
{2}^{-n_{j}} < Col_3^{n_{j_1}}\left( c_{\alpha_i}\right)  < {2}^{-n_{j_1}} + 2\cdot 3^{-k_j-q} . %\cdot {2}^{-n_j}.
\end{equation}
Then,
\begin{equation}\label{key22}
\frac{3^{k_j}}{2^{n_j }} < Col_4^{n_j}\left( c_{\alpha_i}\right) <\frac{3^{k_j}}{2^{n_j }} + \frac{2}{3^{q}}. %\cdot {2}^{-n_j}.
\end{equation}

\vskip 5pt
\noindent
 If $ 1\le \psi_0<2 $, put $q=p+1$. By \eqref{key46} and \eqref{key22},  since $ Col_4^{n_j}\left( c_{\alpha_i}\right) =  c_{\alpha_i+n_j} $, for all $i$,
\begin{equation}\label{key24}
\psi -\frac{1}{3}  < 3^p c_{\alpha_i+n_j}<  \psi +1 .
\end{equation}
Then, for all $i$,
\begin{equation}\label{key25}
 \left[ 3^{p-q_{\alpha_i+n_j}} a_{\alpha_i+n_j} \right] = \left[ 3^p c_{\alpha_i+n_j} \right] = \psi .
\end{equation}

\vskip 5pt
\noindent
If $ 2\le \psi_0<3 $, note that 
\begin{equation}\label{key28}
1< \frac{3^{k_j}}{2^{n_j }} < 2 %\cdot {2}^{-n_j}.
\end{equation}
Choose $q \ge p+2$, such that
\begin{equation}\label{key29}
1<\frac{3^{k_j}}{2^{n_j }} < c_{\alpha_i+n_j} <\frac{3^{k_j}}{2^{n_j }} + \frac{2}{3^{q}} <2 . %\cdot {2}^{-n_j}.
\end{equation}
%By \eqref{key29},  
Since $\left[  c_{\alpha_i+n_j} \right] =1$, then $\left[  Col_3\left( c_{\alpha_i+n_j}\right)  \right] =0$ and
\begin{equation}\label{key100}
 c_{\alpha_i+n_j+1} = Col_4\left(  c_{\alpha_i+n_j} \right) =  3\cdot Col_3\left( c_{\alpha_i+n_j}\right) 
\end{equation}

\begin{equation}\label{key31}
 \frac{3}{2} \cdot c_{\alpha_i+n_j} \le c_{\alpha_i+n_j+1} \le \frac{3}{2} \cdot c_{\alpha_i+n_j} + \frac{3^{-q_{\alpha_i+n_j}}}{2}.
\end{equation}
Then,
\begin{equation}\label{key32}
\frac{3^{k_j+1}}{2^{n_j+1 }} < c_{\alpha_i+n_j+1} <\frac{3^{k_j+1}}{2^{n_j +1}} + \frac{1}{3^{q-1}} + \frac{3^{-q_{\alpha_i+n_j}}}{2}.
\end{equation}
By \eqref{key48},
\begin{equation}\label{key33}
\psi -\frac{1}{3} < 3^p \cdot c_{\alpha_i+n_j+1} <\psi+ 1,
\end{equation}
(since $q_\alpha \rightarrow\infty$, we can assume that $3^{p+1-q_{\alpha_i+n_j}}<1$).

\vskip 5pt
\noindent
Then, for all $i$,
\begin{equation}\label{key35}
\left[ 3^{p-q_{\alpha_i+n_j+1}} a_{\alpha_i+n_j+1} \right] = \left[ 3^p c_{\alpha_i+n_j+1} \right] = \psi .
\end{equation}

%\section{further coments}
\vskip 5pt
\noindent
An inmediate consequence of Theorem \ref{collatz} is the following corollary.

{\cor 
\begin{equation}\label{extra}
\overline{\left\lbrace c_\alpha\right\rbrace }_{\alpha\in \mathbb{N}} =\left[ 1,3\right] .
\end{equation}}

\vskip 5pt
\noindent
Let's define an order relation for integers  $\phi \ge 1$ and $\theta\ge 1 $ in the following way
\begin{equation}\label{extra2}
\phi \preceq \theta \Longleftrightarrow \phi=\left[ 3^{p-q}\theta\right] ,
\end{equation}
where $ p=\left[ \frac{\ln\phi}{\ln3}\right] $ and  $ q=\left[ \frac{\ln\theta}{\ln3}\right] $. We finish by noticing the following corollary.

{\cor  For any integer $\phi\ge 1$, there is a  subsequence $\left(  a_{\alpha_i}\right)_{i\in\mathbb{N}}  $. increasing with respect to the order $\preceq$ defined above% in \ref{extra2}
, minored by $\phi$. }

\bibliographystyle{amsplain}
\bibliography{Collatz}
\end{document}